\documentclass[11pt,amsfonts]{article}
\usepackage{graphicx}
\usepackage{latexsym}
\usepackage{amssymb}
\usepackage{amsmath}
\usepackage{layout}
\newtheorem{prop}{Proposition}
\newtheorem{lemma}{Lemma}

\newtheorem{theorem}{Theorem}
\newtheorem{remark}{Remark}

\def\real{{\mathord{{\rm I\kern-2.8pt R}}}}        
\def\inte{{\mathord{{\rm I\kern-2.8pt N}}}}

\def\sZZ{{\rm Z\kern-2.8ptem{}Z}}

\def\z{{\mathchoice
  {\sZZ}
  {\sZZ}
  {\rm Z\kern-0.30em{}Z}
  {\rm Z\kern-0.25em{}Z} }}
\def\sQQ{{\kern 0.27em \vrule height1.45ex width0.03em depth0em
          \kern-0.30em \rm Q}}
\def\qu{{\mathchoice
    {\sQQ}
    {\sQQ}
  {\kern 0.225em \vrule height1.05ex width0.025em depth0em \kern-0.25em \rm Q}
  {\kern 0.180em \vrule height0.78ex width0.020em depth0em \kern-0.20em \rm Q}
        }}
\def\sCC{{\kern 0.27em \vrule height1.45ex width0.03em depth0em
          \kern-0.30em \rm C}}
\def\complex{{\mathchoice
    {\sCC}
    {\sCC}
  {\kern 0.225em \vrule height1.05ex width0.025em depth0em \kern-0.25em \rm C}
  {\kern 0.180em \vrule height0.78ex width0.020em depth0em \kern-0.20em \rm C}
        }}

\newcommand{\pr}{\mathbb{P}}

\newcommand{\ba}{\begin{array}}
\newcommand{\ea}{\end{array}}
\newcommand{\be}{\begin{equation}}
\newcommand{\ee}{\end{equation}}
\newcommand{\bea}{\begin{eqnarray}}
\newcommand{\eea}{\end{eqnarray}}
\newcommand{\beaa}{\begin{eqnarray*}}
\newcommand{\eeaa}{\end{eqnarray*}}

\def\b{\beta}

\def\z{\zeta}

\font\tenmath=msbm10 \font\sevenmath=msbm7 \font\fivemath=msbm5
\newfam\mathfam \textfont\mathfam=\tenmath
\scriptfont\mathfam=\sevenmath \scriptscriptfont\mathfam=\fivemath

\def \b{\noindent}

\def \={{\buildrel {\rm (law)} \over =}}

\def\qed{ \hfill \vrule width.25cm height.25cm depth0cm\smallskip}

\newcommand{\basa}{\begin{assumption}}
\newcommand{\easa}{\end{assumption}}

\newcommand{\bas}{\begin{assum}}
\newcommand{\eas}{\end{assum}}

\newcommand{\ignore}[1]{}
\begin{document}

\renewcommand{\thefootnote}{\fnsymbol{footnote}}

\title{Maximum likelihood estimators and random walks in long memory models }

\vskip1cm

\author{Karine Bertin $^{1}\quad$ Soledad Torres  $^{1}\quad$ Ciprian A. Tudor $^{2}\vspace*{0.1in}$ \footnote{New affiliation since September 2009: Laboratoire Paul Painlev\'e, Universit\'e de Lille 1,
 F-59655 Villeneuve d'Ascq, France, email: tudor@math.univ-lille1.fr. }
\\$^{1}$  Departamento de Estad\'istica, CIMFAV  Universidad de Valpara\'iso,\\ Casilla
 123-V, 4059 Valparaiso, Chile.
\\soledad.torres@uv.cl\vspace*{0.1in} karine.bertin@uv.cl\vspace*{0.1in}\\$^{2}$SAMOS/MATISSE,
Centre d'Economie de La Sorbonne,\\ Universit\'e de
Panth\'eon-Sorbonne Paris 1,\\90, rue de Tolbiac, 75634 Paris
Cedex 13, France.\\tudor@univ-paris1.fr\vspace*{0.1in}} \maketitle

\begin{abstract}
We consider statistical models driven by Gaussian and non-Gaussian
self-similar processes with long memory and we construct maximum
likelihood estimators (MLE) for the drift parameter. Our approach is
based in the non-Gaussian case on the approximation by random walks
of the driving noise. We study the asymptotic behavior of the
estimators and we give some numerical simulations to illustrate our
results.
\end{abstract}

 \vskip0.5cm

{\bf  2000 AMS Classification Numbers: } 60G18, 62M99.

 \vskip0.3cm


\vskip0.3cm


\section{Introduction}

The self-similarity property for a stochastic process means that scaling of time is equivalent to an appropriate scaling of
space. That is, a process $(Y_{t})_{t\geq 0}$ is selfsimilar of order $H>0$ if for all $c>0$ the processes $(Y_{ct}) _{t \geq
0} $ and $(c^{H} Y _{t} )_{t\geq 0}$  have the same finite dimensional distributions. This property is crucial in applications
such as network traffic analysis, mathematical finance, astrophysics, hydrology  or image processing. We refer to the
monographs \cite{Beran}, \cite{EM} or \cite{ST} for complete expositions on theoretical and practical aspects of self-similar
stochastic  processes.

The most popular self-similar process is the fractional Brownian motion (fBm).  Its practical applications are notorious.  This
process  is defined as a centered Gaussian process $(B^{H}_{t}) _{t\geq 0}$ with covariance function
\begin{equation*}
R^{H}(t,s):=\mathbb{E} (B^{H}_{t}B^{H} _{s}) =\frac{1}{2} \left( t^{2H } + s^{2H} -\vert t-s \vert ^{2H}\right), \hskip0.5cm t,s \geq 0.
\end{equation*}
It can be also  defined as the only Gaussian self-similar process
with stationary increments. Recently, this stochastic process has
been widely studied from the stochastic calculus point of view as
well as from the  statistical analysis point of view. Various types
of stochastic integrals with respect to it have been introduced  and
several  types of stochastic differential equations driven by fBm
have been considered  (see e.g. \cite{N}, Section 5). Another
example of a self-similar process still with long memory (but
non-Gaussian) is the so-called Rosenblatt process  which appears as
limit in limit theorems for stationary sequences with a certain
correlation function (see \cite{DM}, \cite{Ta1}). Although it
received a less important attention than the fractional Brownian
motion, this process is still of interest in practical applications
because of its self-similarity, stationarity of increments and
long-range dependence. Actually the numerous uses of the fractional
Brownian motion in practice (hydrology, telecommunications) are due
to these properties; one prefers in general fBm before other
processes because it is a Gaussian process and the calculus for it
is easier; but in concrete situations when the Gaussian hypothesis
is not plausible for the model, the Rosenblatt process may be an
interesting alternative model. We mention also the work \cite{Taq3}
for examples of the utilisation of non-Gaussian self-similar processes in
practice.

The stochastic analysis of the fractional Brownian motion naturally led to the
statistical inference for  diffusion processes with fBm as driving
noise. We will study in this paper  the problem of the estimation of the drift
parameter. Assume that we have the model
$$dX_t = \theta b(X_t) dt + dB^{H}_t, \hskip0.5cm t\in [0,T]$$
where $(B^{H}_{t})_{t\in [0,T]}$ is a fractional Brownian motion
with Hurst index $H\in (0, 1)$, $b$ is a deterministic function
satisfying some regularity conditions and  the
parameter $\theta \in \mathbb{R}$ has to be estimated. Such
questions have been recently  treated in several papers  (see
\cite{KB}, \cite{TV}  or \cite{SoTu}): in general the techniques used to construct
maximum likelihood estimators (MLE) for the drift parameter
$\theta$ are based on Girsanov transforms for fractional Brownian
motion and depend on  the properties of the deterministic
fractional operators related to the fBm. Generally speaking, the authors of these papers assume that the whole trajectory of the process is continuously  observed. Another possibility is to
use Euler-type approximations for the solution of the above
equation and to construct a MLE estimator based on the density of
the observations given ''the past'' (as in e.g. \cite{Rao},
Section 3.4, for the case of stochastic equations driven by the Brownian motion).

In this work our purpose is to make a first step in the direction of statistical  inference for diffusion processes with
self-similar, long memory and non-Gaussian driving noise. As far as we know, there are not many result on statistical
inference for stochastic differential equations driven by non-Gaussian processes which in addition are not semimartingales.
The basic example of a such process is the Rosenblatt process. We consider here the simple model
\begin{equation*}
X_{t}= at + Z^{H}_{t}, \hskip0.5cm
\end{equation*}
where $(Z^{H}_{t})_{t\in [0,T]}$ is a Rosenblatt process  with known
self-similarity index $H\in (\frac{1}{2}, 1)$ (see Sections
\ref{walkrosen} and Appendix for the definition) and $a\in
\mathbb{R}$ is the parameter to be estimated. We mention that, since
this process is not a semimartingale, it is not Gaussian and its
density function is not explicitly known,  the techniques considered in
the Gaussian case cannot be applied here. We therefore use a
different approach: we consider an approximated model in which we
replace the noise $Z^{H}$ by a two-dimensional disturbed random walk
$Z^{H,n}$ that, from a result is \cite{ToTu}, converges weakly in
the Skorohod topology to $Z^{H}$as $n\to \infty$. Note that this
approximated model still keeps the main properties of the original
model since the noise is asymptotically self-similar and it exhibits
long range dependence. We then construct a MLE estimator (called
sometimes in the literature, see e.g. \cite{Rao} "pseudo-MLE
estimator") using an Euler scheme method  and we prove that this
estimator is  consistent. Although we have not martingales in the
model, this construction involving random walks allows to use
martingale arguments to obtain the asymptotic behavior of the
estimators. Of course, this does not solve the problem of estimating
$a$ in the standard model defined above but we think that our
approach represents a step into the direction of developing models
driven by non-semimartingales and non-Gaussian noises.

Our paper is organized as follows. In Section \ref{prelim} we recall
some facts on the pseudo MLE estimators for the drift parameter in
models driven by the standard Wiener process and by the fBm.  We
construct, in each model, estimators for the drift parameter and we
prove their  strong consistency (in the almost sure sense) or their
$L^ {2}$ consistency under the condition  $\alpha
>1$ where $N^{\alpha}$ is the number of observations at our disposal
and the step of the Euler scheme is $\frac{1}{N}$. This condition
extends the usual hypothesis in the standard Wiener case (see
\ref{c1}, see also \cite{Rao}, paragraph 3.4). Section
\ref{walkrosen} is devoted to the study of the situation when the
noise is the approximated Rosenblatt process; we construct again the
estimator through an inductive method and we study its asymptotic
behavior. The strong consistency is obtained under similar
assumptions as in the Gaussian case. Section \ref{simu} contains
some numerical simulations and in the Appendix we recall the
stochastic integral representations for the fBm and for the
Rosenblatt process.


\section{Preliminaries}\label{prelim}
Let us start by recalling some known facts on maximum likelihood
estimation in simple standard cases. Let $(W_{t}) _{t\in [0,T]}$
be a Wiener process on a classical Wiener space $(\Omega,
{\cal{F}}, P)$ and let us consider the following simple model
\begin{equation}
\label{modW} Y_{t}= at + W_{t}, \hskip0.5cm t\in [0,T]
\end{equation}
with $T>0$ and assume that the parameter $a\in \mathbb{R}$ has to
be estimated. One can for example use the Euler type
discretization of (\ref{modW})
$$Y_{t_{j+1}}^{(n)}:= Y_{t_{j+1}} = Y_{t_j} + a\Delta t + W_{t_{j+1}}-W_{t_j},
\hskip0.5cm j=0, \ldots , N-1, $$ with $Y_{t_0}=Y_{0}=0$ and
$\Delta t= t_{j+1}-t_j$ the step size of  the partition. In the
following, we denote $Y_{t_j}=Y_j$. In the following $f_Z$ denotes
the density of the random variable $Z$.

The conditional density of $Y_{j+1}$ with respect to $Y_{1},
\ldots Y_{j}$ is, by the Markov property,  the same as the
conditional density of $Y_{j+1}$ with respect to $Y_{j}$ and since
$W_{t_{j+1}}-W_{t_j}$ has the normal law $N(0, \Delta t)$, this
density can be expressed by
$$f_{Y_{j+1}/Y_{j}}(y_{j+1}/y_j) = \frac{1}{\sqrt{2\pi (\Delta t)}}
\exp \left(  -\frac{1}{2} \frac{ (y_{j+1} -y_{j} -a\Delta t )
^{2}}{\Delta t } \right).
$$
We easily obtain the likelihood function of the observations
$Y_{1}, \ldots , Y_{N}$

\begin{equation*}
L(a, y_{1}, \ldots , y_{N}) =f_{Y_{1}}(y_1) \prod _{j=1}^{N-1}
f_{Y_{j+1}/Y_{j}}(y_{j+1}/y_j) = \frac{1}{\left( 2\pi (\Delta
t)\right) ^{N/2} } \exp \left( -\frac{1}{2} \sum
_{j=0}^{N-1}\frac{ (y_{j+1} -y_{j} -a\Delta t ) ^{2}}{\Delta t }
\right),
\end{equation*}
and this gives a maximum likelihood estimation of the form
\begin{equation*}
\hat{a}_N= \frac{1}{N\Delta t} \sum _{j=0}^{N-1} (Y_{j+1}-Y_{j}),
\end{equation*}
and the difference $\hat{a}_N-a$ can be written as
\begin{equation*}
\label{a-aW} \hat{a}_N-a= \frac{1}{N\Delta t} \sum_{j=0}^{N-1}
(W_{t_{j+1}}-W_{t_j}).
\end{equation*}
Then
$$\mathbb{E}\left| \hat{a}_N-a\right| ^{2} = \frac{1}{N\Delta t}, $$
and this converges to zero  (that is, the estimator is
$L^{2}$-consistent) if and only if \begin{equation} \label{c1}
 N\Delta t\to \infty, \hskip0.5cm N\to \infty .
\end{equation}
Note that the partition $t_j=\frac{j}{N}$ with $j=0, \ldots N$
does not satisfy (\ref{c1}).

\begin{remark} Under condition (\ref{c1}) we get, by the strong
law of large numbers, the almost sure convergence of the estimator
$\hat{a}_N$ to the parameter $a$.
\end{remark}

\begin{remark}
We need in conclusion to consider an interval between observation
of the order  $\Delta t= \frac{1}{N^{\alpha }} $ with $\alpha <1$
to have (\ref{c1}). Equivalently, if we dispose on $N^{\alpha } $
observations with $\alpha >1$, i.e. $T> N^{\alpha -1}$, and the
interval $\Delta t $ is of order $\frac{1}{N}$, condition
(\ref{c1}) still holds. Using this fact, we will denote in the
sequel by $N^{\alpha} $ the number of observations and we will use
discretization of order $\frac{1}{N}$ of the model.
\end{remark}

\vskip0.5cm

Next, let us take a look to the situation when the Brownian motion
$W$ is replaced by a fractional Brownian motion $B^{H}$ with Hurst
parameter $H\in (0, 1 )$. As far as we know, this situation has not
been considered in the literature. The references \cite{KB} or
\cite{TV} uses a different approach, based on the observation of the
whole trajectory of the process $Y_{t}$ below. The model is now
\begin{equation}
\label{modB} Y_{t}= at+ B^{H}_{t}, \hskip0.5cm t\in [0,T],
\end{equation}
and as before we aim to estimate  the drift parameter $a$ by
assuming that $H$ is known and on the basis on discrete
observations $Y_{1}, \ldots , Y_{N^{\alpha }}$ (the condition on
$\alpha $ will be  clarified later). We use the same Euler type
method  with  $t_j=\frac{j}{N}$ and we denote $Y_{t_j}=Y_j$.

We can easily found the following expression for the observations
$Y_{j}, j=1,\ldots , N^{\alpha}$,
\begin{equation*}
Y_{j}= j\frac{a}{N}+ B^{H} _{\frac{j}{N}}.
\end{equation*}
We need to compute the density of the vector $(Y_{1}, \ldots,
Y_{N^{\alpha}})$. Since the covariance matrix of this vector is
given by $\Gamma= (\Gamma _{i,j})_{i,j=1, \ldots , N^{\alpha }} $
with
$$\Gamma _{i,j}= Cov \left( B^{H}_{\frac{i}{N} }B^{H}_{\frac{j}{N} }\right),$$
the density of $(Y_{1}, \ldots, Y_{N^{\alpha}})$ will be given by
\begin{equation*}
(2\pi ) ^{-\frac{N^{\alpha }}{2}} \frac{1}{\sqrt{ \det \Gamma
}}\exp \left( -\frac{1}{2} \left(y_{1}-\frac{a}{N}, \ldots ,
y_{N^{\alpha }}-N^{\alpha }\frac{a}{N}\right) ^{t} \Gamma^{-1}
\left(y_{1}-\frac{a}{N}, \ldots , y_{N^{\alpha }}-N^{\alpha
}\frac{a}{N}\right) \right)
\end{equation*}
and by maximizing the above expression with respect to the variable $a$ we obtain the following MLE estimator
\begin{equation}
\label{MLEBH} \hat{a}_{N}= N\frac{\sum_{i,j=1}^{N^{\alpha }}
j\Gamma _{i,j} ^{-1} Y_{i}  }{\sum_{i,j=1}^{N^{\alpha }}ij \Gamma
_{i,j} ^{-1}}
\end{equation}
and then
\begin{equation}
\label{an-aBH} \hat{a}_{N} -a = N\frac{\sum_{i,j=1}^{N^{\alpha }}
j\Gamma _{i,j} ^{-1} B^{H}_{\frac{i}{N}} }{\sum_{i,j=1}^{N^{\alpha
}}ij \Gamma _{i,j} ^{-1}},
\end{equation}
where the $\Gamma _{i,j} ^{-1}$ are the coordinates of the matrix
$\Gamma^{-1}$.

\noindent Thus
\begin{eqnarray*}
\mathbb{E}\left| \hat{a}_N-a\right| ^{2} &=&
N^{2}\frac{\sum_{i,j,k,l=1}^{N^{\alpha}}jl\Gamma ^{-1}_{i,j}
\Gamma ^{-1}_{k,l}
E\left( B^{H} _{\frac{i}{N}}B^{H} _{\frac{k}{N}} \right) }{\left( \sum_{i,j=1}^{N^{\alpha}}ij\Gamma ^{-1}_{i,j}\right) ^{2}}\\
&=&N^{2}\frac{\sum_{i,j,k,l=1}^{N^{\alpha}}jl\Gamma ^{-1}_{i,j}
\Gamma ^{-1}_{k,l} \Gamma _{i,k} }{\left(
\sum_{i,j=1}^{N^{\alpha}}ij\Gamma ^{-1}_{i,j}\right) ^{2}}.
\end{eqnarray*}
Note that
\begin{equation*}
\sum_{i,j,k,l=1}^{N^{\alpha}}jl\Gamma ^{-1}_{i,j} \Gamma
^{-1}_{k,l} E\left( B^{H} _{\frac{i}{N}}B^{H} _{\frac{k}{N}}
\right) =\sum_{j,k,l=1}^{N^{\alpha }}jl\Gamma _{k,l} ^{-1}  \left(
\sum_{i=1}^{N^{\alpha }}\Gamma _{i,j}^{-1}\Gamma _{i,k}
\right)=\sum_{j,k,l=1}^{N^{\alpha }}jl\Gamma _{k,l} ^{-1}  \delta
_{jk} =  \sum_{j,l=1}^{N^{\alpha}}jl\Gamma ^{-1}_{j,l}
\end{equation*}
and consequently
\begin{equation*}
E\left| \hat{a}_{N} -a \right| ^{2} = N^{2} \frac{1}{
\sum_{i,j=1}^{N^{\alpha}}ij\Gamma ^{-1}_{ij}  }= N^{2-2H}
\frac{1}{ \sum_{i,j=1}^{N^{\alpha}}ijm ^{-1}_{ij}  }
\end{equation*}
where $m^{-1}_{ij}$ are the coefficients of the matrix $M^{-1}$  with $M= (m_{ij})_{i,j=1, \ldots , N^{\alpha }}$, $m_{ij} =
\frac{1}{2} (i^{2H} + j^{2H}-\vert i-j\vert ^{2H} ) $. Let $x$ be the vector $(1,2, \ldots , N^{\alpha })$ in
$\mathbb{R}^{N^{\alpha }}$. We use the inequality
\begin{equation*}
x^{t} M^{-1} x \geq \frac{\Vert x\Vert ^{2}_{2} }{\lambda }
\end{equation*}
where $\lambda $ is the largest eigenvalue of the matrix $M$. Then we will have
\begin{equation*}
\mathbb{E}\left| \hat{a}_{N} -a \right| ^{2}  \leq N^{2-2H} \frac{\lambda }{\Vert x\Vert _{2}^{2}}.
\end{equation*}
Since $1^{2} + 2^{2} +\ldots p^{2}= \frac{p(p+1)(2p+1)}{6} $ we
see that $\Vert x\Vert _{2}^{2}$ behaves as $N^{3\alpha }$ and on
the other hand (by the Gersghorin Circle Theorem (see \cite{GvL},
Theorem 8.1.3, pag. 395)
\begin{equation*}
\lambda \leq \max _{ i=1, ..., N^{\alpha }}\sum_{l=1}^{N^{\alpha
}} \vert m_{il}\vert \leq C N^{\alpha (2H+1)},
\end{equation*}
with $C$ a positive constant. Finally,
\begin{equation}\label{est1}
\mathbb{E}\left| \hat{a}_{N} -a \right| ^{2}  \leq
CN^{2-2H-3\alpha + \alpha (2H-1)}=C N^{(2-2H)(1-\alpha )}
\end{equation}
and this goes to zero if and only if $\alpha >1$.

\vskip0.2cm

Let us summarize the above discussion.

\begin{prop}
Let $(B^{H}_{t})_{t\in [0,T]}$  be a fractional Brownian motion with Hurst parameter $H\in (0,1)$ and let $\alpha >1$.  Then the estimator (\ref{MLEBH}) is $L^{p}$ consistent for any $p\geq 1$.
\end{prop}
{\bf Proof: } Since for every $n$ the random variable $\hat{a}_{n} -a$ is a centered random variable, it holds that, for some positive constant depending on $p$
\begin{equation*}
\mathbb{E}\left| \hat{a}_{N} -a \right| ^{p}  \leq c_{p} \left( \mathbb{E}\left| \hat{a}_{N} -a \right| ^{2} \right) ^{\frac{p}{2}} \leq c_{p}N^{p(1-H)(1-\alpha)}
\end{equation*}
and this converges to zero as zero since $\alpha >1$. \qed

\vskip0.5cm

It is also possible to obtain the almost sure convergence of the estimator to the the true parameter from the the estimate (\ref{est1}).

\begin{prop}
Let $(B^{H}_{t})_{t\in [0,T]}$  be a fractional Brownian motion with Hurst parameter $H\in (0,1)$ and let $\alpha >1$. Then the estimator (\ref{MLEBH}) is strongly consistent, that is, $\hat{a}_{n} -a \to _{n\to \infty}p$ almost surely.
\end{prop}
{\bf Proof: } Using Chebyshev's inequality,
\begin{equation*}
P\left( \vert \hat{a}_{n}-a\vert >N^{-\beta }\right) \leq c_{p}N^{\beta p} N^{p(1-H)(1-\alpha)}.
\end{equation*}
In order to apply the Borel-Cantelli lemma, we need to find a strictly positive $\beta$ such that $\sum_{N\geq 1}N^{\beta p} N^{p(1-H)(1-\alpha)}<\infty.$  One needs $p\beta  +(1-H)(1-\alpha )p <-1$ and this is  possible if and only if  $\alpha >1$. \qed

\vskip0.3cm

\begin{remark}
The fact that the restriction $\alpha >1$ is interesting because it
justifies looking for another type of estimator for which no such
restriction holds. But for the approach using the Euler scheme
discretization, this restriction is somehow expected  because it
appear also also in the Wiener case (see \cite{Rao}).
\end{remark}

\begin{remark}
Let us also comment on the problem of estimation of the diffusion
parameter in the model (\ref{modB}). Assume that the fractional
Brownian motion $B^{H}$ is replaced by $\sigma B^{H}$ in
(\ref{modB}), with $\sigma \in \mathbb{R}$. In this case it is
known that the sequence $$N^{2H-1} \sum_{i=0}^{N-1} \left(
Y_{\frac{i+1}{N}} -Y_{\frac{i}{N}}\right) ^{2} $$ converges (in
$L^{2} $ and almost surely ) to $\sigma ^{2}$. Thus we easily
obtain an estimator for the diffusion parameter by using such
quadratic variations. The above sequence has the same behavior if
we replace the fBm $B^{H}$ by the Rosenblatt process $Z^{H}$
because the Rosenblatt process is also self-similar with
stationary increments and it still  satisfies $\mathbb{E}\left|
Z^{H}_{t} -Z^{H}_{s}\right| ^{2} = \vert t-s\vert ^{2H}$ for every
$s,t\in [0,T]$ (see Section 4 for details). For this reason we
assume throughout this paper that the diffusion coefficient is
equal to $1$.
\end{remark}

\section{MLE and random walks in the non-Gaussian case}\label{walkrosen}

We study in this section a non-Gaussian long-memory model. The
driving process is now a Rosenblatt process with selfsimilarity
order $H\in (\frac{1}{2}, 1)$. This process appears as a limit in
the so-called {\it Non Central Limit Theorem } (see \cite{DM} or
\cite{Ta1}). It can be defined through its representation as
double iterated integral with respect to a standard Wiener process
given by equation (\ref{repZ}) in the Appendix. Among its main
properties, we recall
\begin{description}
\item{$\bullet$ } it exhibits long-range dependence (the
covariance function decays at a power function at zero)

\item{$\bullet$ }it is $H$-selfsimilar in the sense that for any
$c>0$, $(Z^{H}(ct)) =  ^{(d)} (c^{H} Z^{H} (t))$, where $ "
=^{(d)} "$ means equivalence of all finite dimensional
distributions ; moreover, it has stationary increments, that is,
the joint distribution of $(Z^{H}(t+h)-Z^{H}(h), t\in [0,T])$ is
independent of $h> 0$.

\item{$\bullet$ } the covariance function is
$$\mathbb{E}(Z^{H}_{t}Z^{H}_s) = \frac{1}{2}\left( t^{2H} + s^{2H } -\vert
t-s \vert ^{2H} \right), \hskip0.5cm s,t\in [0,T] $$

and consequently, for every $s,t\in [0,T]$
\begin{equation*}
\mathbb{E}\left| Z^{H}_{t} -Z^{H}_s\right| ^{2} = \vert t-s\vert
^{2H}
\end{equation*}

\item{$\bullet$ }  the Rosenblatt  process is H\"older continuous
of order $\delta <H$

\item{$\bullet$ } it is not a Gaussian process; in fact, it can be
written as a double stochastic integral of a two-variable
deterministic function with respect to the Wiener process.

\end{description}

Assume that we want to construct a MLE estimator for the drift
parameter $a$ in the model
$$Y_{t}= at+ Z^{H}_{t}, \hskip0.5cm t\in [0,T]. $$
We first note that the Rosenblatt process has only been defined for
the self-similarity order $H>\frac{1}{2}$ and consequently we will
have in the sequel always $H\in (\frac{1}{2}, 1)$.  The approaches
used previously do not work anymore because the Rosenblatt process
is still not a semimartingale and moreover, in contrast  to the fBm
model, its density function is not explicitly known anymore. The
method based on random walks approximation offers a solution to the
problem of estimating the drift parameter $a$. We will use this
direction; that is, we will replace the process $Z_H$ by its
associated random walk
\begin{equation}
\label{zn} Z^{H,N}_{t}= \sum_{i,k=1; i\not= k}^{[Nt]}N^{2}
\int_{\frac{i-1}{N} }^{\frac{i}{N}}\int_{\frac{k-1}{N}}^
{\frac{k}{N}}F\left( \frac{[Nt]}{N}, u,v\right) dvdu \frac{\xi
_{i}}{\sqrt{N}}\frac{\xi _{k}}{\sqrt{N}}, \hskip0.5cm t\in [0,T]
\end{equation}
where the $\xi_i$ are i.i.d. variables of mean zero and variance one
and the deterministic kernel $F$ is defined in Appendix by
(\ref{defF}). It has been proved in \cite{ToTu} that the random walk
(\ref{zn}) converges weakly in the Skorohod topology to the
Rosenblatt process $Z^H$.\\

We consider the following discretization of the Rosenblatt process
\begin{equation*}(Z^{N,H} _{\frac{j}{N}}), \quad j=0, \ldots N^{\alpha},
\end{equation*}
where for $j\neq 1$ $Z^{H,N} _{\frac{j}{N}}$ is given by
(\ref{zn}) and we set $Z^{N}_{\frac{1}{N}}=\xi_1/N^{H}$.
 With this slight modification, the process $(Z^{H,N}
_{\frac{j}{N}})$ still converges weakly in the Skorohod topology to
the Rosenblatt process $Z^H$. We will assume as above that the
variables $\xi _i $ follows a standard normal law $N(0,1)$.

Concretely, we want to  estimate the drift parameter $a$ on the
basis of the observations
\begin{equation*}
 Y_{t_{j+1}} = Y_{t_j}+ a (t_{j+1}-t_j) + \left(
Z^{H,N}_{t_{j+1}}-Z^{H,N}_{t_j}\right)
\end{equation*}
where $t_j=\frac{j}{N}$, $j=0,\ldots , N^{\alpha } $ and
$Y_{0}=0$. We will assume again that we have at our disposal a
number $N^{\alpha }$ of observations and we use a discretization
of order $\frac{1}{N}$ of the model. Denoting $Y_{j}:= Y_{t_j}$,
we can write \begin{equation*} Y_{j+1}= Y_{j} + \frac{a}{N}
+\left( Z^{H,N}_{\frac{j+1}{N}}-Z^{H,N}_{\frac{j}{N}}\right),
\hskip0.5cm j=0, \ldots, N^{\alpha }-1.
\end{equation*}
Now, we have
\begin{eqnarray*}
Z^{H,N}_{\frac{j+1}{N}}-Z^{H,N}_{\frac{j}{N}} &=&  f_{j}(\xi _{1},
\ldots , \xi_{j}) + g_{j}(\xi_1,\ldots,\xi_j)\xi _{j+1}
\end{eqnarray*}
where $f_0=0$, $f_1=f_1(\xi_1)=-\xi_1/N^{H}$ and for $j\geq 2$
\begin{equation*}
f_j=f_{j}(\xi_{1}, \ldots , \xi_{j})= N\sum _{i,k=1; i\not= k}^{j}
\left( \int _{\frac{i-1}{N}}^{\frac{i}{N}}\int
_{\frac{k-1}{N}}^{\frac{k}{N}}\left( F\left(
\frac{j+1}{N},u,v\right) -F\left( \frac{j}{N},u,v\right)\right)
dvdu\right) \xi_{i}\xi_{k},
\end{equation*}
and $g_0=1/N^{H}$ for $j\ge 1$
\begin{equation*}
g_j=g_{j} (\xi_{1}, \ldots , \xi_{j}) = 2 N \sum_{i=1}^{j} \left(
\int _{\frac{i-1}{N}}^{\frac{i}{N}}\int
_{\frac{j}{N}}^{\frac{j+1}{N}}F\left( \frac{j+1}{N},
u,v\right)dvdu \right) \xi_{i}.
\end{equation*}
Finally we have the model
\begin{equation}\label{modR}
Y_{j+1}= Y_{j} + \frac{a}{N}+f_{j} + g_{j}\xi _{j+1}.
\end{equation}

In the following, we assume that $(\xi_1,\ldots,\xi_n,\ldots)\in
B$ where
$$B=\cap_{j\ge 1}\left\{g_{j}(\xi_1,\ldots,\xi_{j})\neq 0\right\}.$$
The event $B$ satisfy $\pr(B)=1$.

\begin{remark}
Note that on the event $B$, conditioning with respect to the
$\xi_{1},\cdots,\xi_{j}$ is the same as conditioning with respect
to the $Y_{1},\cdots,Y_{j}$. In fact, we have
\begin{eqnarray}
Y_1 &=& a/N +  \xi_1/N^H 
\nonumber \\
Y_j &=& Y_{j-1}+a/N +  f_{-1}(\xi_{1},\ldots , \xi_{j-1})+ \xi_{j}
 g_{j-1}(\xi_{1},\ldots , \xi_{j-1}), \quad j\ge 2
\nonumber.
\end{eqnarray}
Since on $B$, for all $j\ge 2$,
$g_{j-1}(\xi_1,\ldots,\xi_{j-1})\neq 0$, the two $\sigma$-algebra
$\sigma(\xi_1,\ldots,\xi_{j-1})$ and $\sigma(Y_1,\ldots,Y_{j-1})$
satisfy
$$\sigma(\xi_1,\ldots,\xi_{j-1})\cap B=\sigma(Y_1,\ldots,Y_{j-1})\cap B.$$
\end{remark}
\ \\

Then, given $\xi_{1}, \ldots , \xi_{j}$ the random variable
$Y_{j+1} $ is conditionally Gaussian and the conditional density
of $Y_{j+1}$ given  $Y_{1}, \ldots , Y_{j}$ can be written as
\begin{equation*}
 f_{Y_{j+1} / Y_{1}, \ldots , Y_{j}}(y_{j+1} / y_{1}, \ldots ,
y_{j}) = \frac{1}{ \sqrt{2\pi g_{j}^{2}}}\exp \left(
-\frac{1}{2}\frac{ (y_{j+1}-y_{j}-a/N
 -f_{j} )^{2}} {g_{j}^{2}}\right).
\end{equation*}
The likelihood function can be expressed as \begin{eqnarray*}
 L(a, y_{1}, \ldots,y_{N^{\alpha }})&=& f_{Y_{1}}(y_1)f_{Y_{2}/Y_{1}}(y_2/y_1)
  \ldots f_{Y_{N^{\alpha }} / Y_{1}, \ldots, Y_{ N^{\alpha }-1} }(y_{N^{\alpha }} / y_{1}, \ldots, y_{ N^{\alpha }-1})\\
 &=& \prod_{j=0}^{N^{\alpha }-1}\frac{1}{
\sqrt{2\pi g_{j}^{2}}}\exp \left(  -\frac{1}{2}\frac{
(y_{j+1}-y_{j}-a/N
 -f_{j} )^2} {g_{j}^{2}}\right).
\end{eqnarray*}

\vskip0.5cm

By standard calculations, we will obtain
\begin{equation}
\label{hatZ} \hat{a}_N=\frac{  N \sum_{j=0}^{N^{\alpha } -1}
\frac{ (Y_{j+1} -Y_{j} -f_{j}-a/N)}{g_{j}^{2}} }{ \sum
_{j=0}^{N^{\alpha } -1} \frac{ 1}{g_{j}^{2}} }
\end{equation}
and since $Y_{j+1} -Y_{j} -f_{j}= \frac{a}{N}+ g_{j}\xi_{j+1} $ we
obtain
\begin{equation*}
\label{a-aZ} \hat{a}_N-a = \frac{N\sum_{j=0}^{N^{\alpha } -1}
\frac{ \xi _{j+1}}{g_{j}} }{ \sum _{j=0}^{N^{\alpha } -1} \frac{
1}{g_{j}^{2}} }.
\end{equation*}

Let us  comment on the above expression. Firstly, note that since
$\sum _{j=0}^{N^{\alpha } -1} \frac{ 1}{g_{j}^{2}}$ is not
deterministic anymore, the square mean of the difference
$\hat{a}_N-a$ cannot be directly computed. Secondly, if we denote
again by
\begin{equation*}
A_{M} =\sum_{j=0}^{M-1} \frac{ \xi _{j+1}}{g_{j}}
\end{equation*}
this discrete process still satisfies
$$\mathbb{E}\left( A_{M+1} / {\cal{F}}_{M} \right) = A_{M}, \hskip0.5cm \forall
M\ge 1
$$
where $\mathcal{F}_M$ is the $\sigma$-algebra generated by
$\xi_1,\ldots,\xi_M$. However we cannot speak about square
integrable martingales brackets because $A_M ^{2}$ is not
integrable (recall that the bracket is defined in general for
square integrable martingales). In fact, $\mathbb{E}(A_{M}^2)
=\sum_{j=0}^{M}\mathbb{E}\left(\frac{ 1}{g_{j}^{2}}\right) $ and
this is not finite in general because $g_{j}$ is a normal random
variable. We also mention that, in contrast to the Gaussian case,
the expectation of the estimator is not easily calculable anymore
to decide if (\ref{hatZ}) is unbiased. On the other hand, from the
numerical simulation it seems that the estimator is biased.

\vskip0.3cm

Nevertheless, martingale type methods can be employed to deal with
the estimator (\ref{hatZ}).

We use again the notation
$$\langle A\rangle _{M}=\sum_{j=0}^{M-1}\frac{
1}{g_{j}^{2}} .$$

The following lemma is crucial.

\begin{lemma}
\label{Linfty2} Assume that $\alpha >2-2H $ and let us denote by
$$T_{N}:= \frac{1}{N^{2}} \sum_{j=0}^{N^{\alpha
}-1}\frac{1}{g_{j}^{2} }.
$$ Then $T_{N} \xrightarrow{N\to \infty } \infty $ almost surely.
Denote
$$U_N=\frac{1}{N^{2}} \left( \sum_{j=0}^{N^{\alpha }-1}\frac{1}{g_{j}^{2} }\right)^{1-\gamma }.$$
Then there exists $0< \gamma <1$, such that $U_{N}
\xrightarrow{N\to \infty } \infty $ almost surely.

\end{lemma}
{\bf Proof: } Let prove the convergence for $T_N$. We will use a
Borel-Cantelli argument. To this end, we will show that
\begin{equation}
\label{toprove} \sum_{N\geq 1} \pr\left( T_{N} \leq N^{\delta }
\right)  <\infty
\end{equation}
for some $\delta >0$.

Fix $0<\delta <\alpha -(2-2H)$.   We have
\begin{eqnarray*}
 \pr\left( T_{N} \leq N^{\delta } \right) &=& \pr\left(\sum_{j=0}^{N^{\alpha }-1}\frac{1}{g_{j}^{2} }\leq N^{2+\delta } \right) \\
 &=& \pr\left( \frac{1} {\sum_{j=0}^{N^{\alpha }-1}\frac{1}{g_{j}^{2} } }\geq N^{-2-\delta} \right) .
\end{eqnarray*}
We now choose a $p$ integer large enough such that $\frac{1}{p} <
\alpha -(2-2H)-\delta$ and we apply the Markov inequality. We can
bound $ \pr\left( T_{N} \leq N^{\delta } \right)$ by
$$ \pr\left( T_{N} \leq N^{\delta } \right)\leq N^{(2+\delta )p}\mathbb{E}\left| \frac{1} {  \sum_{j=0}^{N^{\alpha }-1}\frac{1}{g_{j}^{2} } }  \right| ^{p}
$$
and using the inequality "harmonic mean is less than the arithmetic mean" we obtain
$$ \pr\left( T_{N} \leq N^{\delta } \right)\leq N^{(2+\delta )p} N^{-2\alpha p} \mathbb{E}\left(  \sum_{j=0}^{N^{\alpha } -1} g_{j}^{2}
\right) ^{p}
$$
Note that (the first inequality below can be argued similarly as in
\cite{Sot})
$$\frac{1}{N^{2H}} \geq\mathbb{E}\left|  Z^{H,N}_{\frac{j+1}{N}}- Z^{H,N}_{\frac{j}{N}} \right| ^{2}
= \mathbb{E}\left(f_{j} ^{2}\right) +
\mathbb{E}\left(g_{j}^{2}\right),
$$
it holds that
\begin{equation*}
\label{imp2} \mathbb{E}\left(g_{j}^{2}\right):=c_{j}\leq
\frac{1}{N^{2H}}.
\end{equation*}

Since $g_{j}= \sqrt{c_{j}} X_{j}$ where $X_{j} \sim  N(0,1)$, by (\ref{imp2})

\begin{eqnarray*}
\pr\left( T_{N} \leq N^{\delta } \right)&\leq & N^{(2+\delta -2\alpha -2H )p}\mathbb{E} \left(  \sum_{j=0}^{N^{\alpha}-1} X_{j}^{2} \right) ^{p} \\
&\leq &N^{(2+\delta -2\alpha -2H )p} N^{\alpha p} = N^{p( 2+\delta
-\alpha -2H)} \end{eqnarray*} and thus relation (\ref{toprove}) is
valid.

The convergence for $U_N$ can be obtained in a similar way with
$0<\gamma<1$ such that $\alpha(2-\gamma)-2+2H\gamma>0$,
$0<\delta<\alpha(2-\gamma)-2+2H\gamma $ and $p$ such that
$1/p<\alpha(2-\gamma) -(2-2H\alpha)-\delta$.
 \qed

We use the notation

$$V_{M}:= \frac{A_{M}^{2}}{\langle A \rangle _{M} ^{1+ \gamma }}, \hskip0.5cm M\geq 1 $$
and
$$B_{M} := \frac{  \langle A\rangle _{M+1} - \langle A\rangle _{M}}{\langle A\rangle _{M+1} ^{1+\gamma }}.$$
Note that $V_M$ and $B_M$ are $\mathcal{F}_M$ adapted.

We recall
the Robbins-Siegmund criterium for the almost sure convergence which will play an important role in the sequel
(see e.g. \cite{duflo}, page 18): let $(V_{N})_{N}, (B_{N})_{N}$
be ${\cal{F}}_{N} $ adapted, positive sequences such that
$$\mathbb{E}\left(  V_{N+1}/ {\cal{F}}_{N}\right) \leq V_{N} + B_{N} \hskip0.5cm \mbox{ a.s. }.
$$
Then the sequence of random variables $(V_{N})_{N}$ converges as
$N\to \infty$ to a random variable
 $V_{\infty }$ almost surely on the set $\{ \sum_{N\geq 1} B_{N} <\infty \}$.
\begin{lemma}
\label{last} The sequence $V_{M}$ converges to a random variable
almost surely when  $M\to \infty$.
\end{lemma}
{\bf Proof: }  We make use of the Robbins-Siegmund criterium. It holds that
\begin{eqnarray*}
&&\mathbb{E}\left( V_{M+1} /{\cal{F}}_{M} \right) = \mathbb{E}
\left( \frac{A_{M+1} ^{2} }{\langle A \rangle _{M+1} ^{1+ \gamma
}}/{\cal{F}}_{M}\right) \\
&& \leq \frac{1}{\langle A \rangle _{M+1} ^{1+ \gamma
}}\mathbb{E}\left( A_{M+1}^2 /{\cal{F}}_{M}\right)   =
\frac{1}{\langle A \rangle _{M+1} ^{1+ \gamma }}\left(  A_{M}^{2}
+ \langle A\rangle _{M+1} -\langle A\rangle _{M}\right)
\\
&&\leq V_{M} + B_{M}.
\end{eqnarray*}
By Lemma \ref{Linfty2} the sequence $\langle A\rangle _{M}$
converges to $\infty$ and therefore
$$\sum _{N} B_{N^{\alpha }} \leq C + \int_{1}^{\infty} x^{-1-\gamma } ds <\infty.$$
Once can conclude by applying the Robbins-Siegmund criterium. \qed

\vskip0.5cm

We state now the main result of this section.
\begin{theorem}
The estimator (\ref{hatZ}) is strongly consistent.
\end{theorem}
{\bf Proof: }We have $$(\hat{a}_N-a)^{2}= \frac{V_{N^{\alpha}}}
{\frac{1}{N^2}\langle A\rangle _{N^\alpha}^{1-\gamma}}
=\frac{V_{N^{\alpha}}}{U_{N}}$$ and we conclude by Lemmas
\ref{Linfty2} and \ref{last}. \qed

\vskip0.5cm

\b {\bf Comment: } The $L^{1}$ (or $L^{p}$) consistence of the estimator  (\ref{hatZ}) is on open problem. Note that its expression is a fraction whose numerator and denominator are non-integrable involving inverses of Gaussian random variables. The natural approaches to deal with do not work. A first basic idea when is to try to apply the H\"older inequality for the product $FG$ with $F=\sum_{j=0}^{N^{\alpha } -1}
\frac{ \xi _{j+1}}{g_{j}}$ and $G=\sum _{j=0}^{N^{\alpha } -1} \frac{
1}{g_{j}^{2}} $. But when we do this, we are confronted with the moments of the random variable $\vert F\vert $ and it has not moments. Indeed, the second moment is $\mathbb{E}\sum _{j=0}^{N^{\alpha } -1} \frac{
1}{g_{j}^{2}} $ which is obviously infinite.  Even $\mathbb{E}\vert F\vert $ is infinite.
Another way is to use the inequality $\left( \sum_{i=1}^{n} \vert a_i  b_{i} \vert  \right) ^{2}\leq \left( \sum_{i=1}^{n} \vert a_i \vert ^{2} \right) \left( \sum_{i=1}^{n} \vert b_i \vert ^{2} \right)$. In this way we avoid the non-integrability of the random variable $F$ defined above.
We get
\begin{equation*}\sum_{j=0}^{N^{\alpha } -1}
\frac{ \xi _{j+1}}{g_{j}} \leq  \left(\sum_{j=0}^{N^{\alpha } -1}\xi_{j+1}^{2}\right) ^{\frac{1}{2}}\left(\sum _{j=0}^{N^{\alpha } -1} \frac{
1}{g_{j}^{2}} \right) ^{\frac{1}{2}}
\end{equation*}
and then
\begin{equation*}
\vert \hat{a}_N-a \vert \leq N\frac{
  \left(\sum_{j=0}^{N^{\alpha } -1}\xi_{j+1}^{2}\right) ^{\frac{1}{2}}}{\left(\sum _{j=0}^{N^{\alpha } -1} \frac{
1}{g_{j}^{2}} \right) ^{\frac{1}{2}}}
\end{equation*}
and by saying that the harmonic mean is less that the arithmetic mean (we have the feeling that it is rather optimal in this case)
\begin{equation}\label{e1}
\vert \hat{a}_N-a \vert\leq  N^{1-\alpha }\left(\sum_{j=0}^{N^{\alpha } -1}\xi_{j+1}^{2}\right) ^{\frac{1}{2}}\left(\sum_{j=0}^{N^{\alpha } -1}g_{j}^{2} \right) ^{\frac{1}{2}}.
\end{equation}
We tried now  to  apply H\"older for $ \mathbb{E}\vert \hat{a}_N-a
\vert$ and this gives
\begin{equation*}
\mathbb{E}\vert \hat{a}_N-a \vert\leq  N^{1-\alpha }\left(
\mathbb{E}\left(\sum_{j=0}^{N^{\alpha } -1}\xi_{j+1}^{2}\right)
^{\frac{p}{2}}\right)^{1/p}\left(
\mathbb{E}\left(\sum_{j=0}^{N^{\alpha } -1}g_{j}^{2} \right)
^{\frac{q}{2}}\right) ^{1/q}
\end{equation*}
But now $\sum_{j=0}^{N^{\alpha } -1}\xi_{j+1}^{2}$ is a chi-square
variable with $N^{\alpha }$ degree of freedom, its $p$ moment
behaves as $N^{\alpha p}$ and it can be also proved that
$\left(\mathbb{E}\left(\sum_{j=0}^{N^{\alpha } -1}g_{j}^{2} \right)
^{\frac{q}{2}}\right) ^{1/q}$ is of order less than $N^{-\alpha }{2}
N^{-H}$ (in the paper it is proved in the case $q=2$) which
unfortunately gives $\sum_{j=0}^{N^{\alpha } -1}\xi_{j+1}^{2}\leq
N^{1-H}\to \infty. $ One needs to compute directly the norm of
$\vert \hat{a}_N-a \vert$ in (\ref{e1}) by using the distribution of
the vector $(g_{1},\ldots , g_{n})$. But  it turns out that this
joint distribution is complicated and not tractable because of the
presence of the kernel of the Rosenblatt process.

\vskip0.5cm


\section{Simulation}\label{simu}
We consider the problem of estimating $a$ from the observations
$Y_1,\ldots,Y_N$ \eqref{modR} driven by the approximated
Rosenblatt process respectively. In all of the cases, we use
$\alpha=2$. We have implemented the estimator $\hat{a}_N$. We have
simulated the observations $Y_1,\ldots,Y_{N^2}$ for different
values of $H: 0.55, 0.75$ and $0.9$ and the values of $a$: 2 and
20. We consider the cases $N=100$ and $N=200$, in others words in
the first case we have $10000$ observations and in the second case
$40000$. For each case, we calculate $100$ estimation $\hat{a}_N$
and we give in the following tables the mean and the standard
deviation of these estimation.

Finally in the model driven by the approximated Rosenblatt process
(Section \ref{walkrosen}), we only construct an estimator that
depends on the observations $Y_j$ and of the $\xi_j$.

The results for $N=100$ are the following.
$$
\begin{array}{|r|r|r|r|}\hline
a=2&H=0.55&   H=0.75&  H=0.9\\\hline \text{mean}& 2.1860& 2.1324&
2.1416\\\hline \text{stand. deviation} & 0.3757 & 0.3702& 0.3393
\\\hline
\end{array}
$$

$$
\begin{array}{|r|r|r|r|}\hline
a=20&H=0.55&   H=0.75&  H=0.9\\\hline \text{mean}& 20.2643&
20.2013& 20.1506\\\hline \text{stand. deviation} & 0.4272 &
0.4982& 0.4507
\\\hline
\end{array}
$$

The results for $N=200$ are the following.

$$
\begin{array}{|r|r|r|r|}\hline
a=2&H=0.55&   H=0.75&  H=0.9\\\hline \text{mean}& 2.0433& 2.1048&
2.0559\\\hline \text{stand. deviation} & 0.2910 & 0.2795& 0.2768
\\\hline
\end{array}
$$

$$
\begin{array}{|r|r|r|r|}\hline
a=20&H=0.55&   H=0.75&  H=0.9\\\hline \text{mean}& 20.1895&
20.1080& 20.0943\\\hline \text{stand. deviation} & 0.1952 &
0.2131& 0.2162
\\\hline
\end{array}
$$

We can observe that the quality of estimation increases as $N$
increases. We obtain these last estimations for the parameter of
discretization $N=200$. Since the computational cost is quite
important, we do not implement the estimator for larger $N$. But
although with $N=200$, the estimator is quite good.

\vspace{0.5cm}

\noindent{\Large \textbf{Appendix: Representation of fBm and
Rosenblatt process as stochastic integral with respect to a Wiener
process}}

\vspace{0.5cm}
 The
fractional Brownian process  $(B^{H}_{t})_{t\in [0,T]}$ with Hurst
parameter $H\in (0,1)$ can be written (see e.g. \cite{N}, Chapter 5
or \cite{N1})
$$B_t^H=\int_0^tK^{H}(t,s)dW_s, \quad t\in[0,T]$$
where $(W_t,t\in[0,T])$ is a standard Wiener process,
\begin{equation}\label{defK}K^H(t,s)=c_Hs^{\frac{1}{2}-H}\int_s^t(u-s)^{H-\frac{3}{2}}u^{H-\frac{1}{2}}du\end{equation}
where $t>s$ and
$c_H=\left(\frac{H(2H-1)}{\beta(2-2H,H-\frac{1}{2})}\right)^{\frac{1}{2}}$
and $\beta(\cdot,\cdot)$ is the beta function. For $t>s$, we have
\begin{equation*}
\label{dK} \frac{\partial K^H}{\partial t} (t,s)= c_{H}\left(
\frac{s}{t} \right) ^{\frac{1}{2}-H} (t-s)^{H-\frac{3}{2}}.
\end{equation*}
An analogous representation for the Rosenblatt process $(Z^{H}_{t})
_{t\in [0,T]}$ is (see \cite{Tu})
\begin{equation}
\label{repZ} Z^{H}_{t}=\int _{0}^{t}\int_{0}^{t} F(t,y_1,y_2)dW_ {y_{1}} dW_{y_{2}}
\end{equation}
where $(W_{t}, t \in [0,T])$ is a Brownian motion,
\begin{equation}\label{defF}F(t,y_1,y_2)= d(H)1_{[0,t]}(y_1)1_{[0,t]}(y_2)\int_{ y_{1} \vee y_{2}
}^{t} \frac{\partial K^{H'} }{\partial u} (u,y_{1} ) \frac{\partial
K^{H'} }{\partial u} (u,y_{2} )du,\end{equation} $ H'=\frac{H+1}{2}$
and $ d(H)= \frac{1}{H+1} \left( \frac{H}{2(2H-1)} \right)
^{-\frac{1}{2}}.$ Actually the representation (\ref{repZ}) should be
understood as a multiple integral of order 2 with respect to the
Wiener process $W$. We refer to \cite{N} for the construction and
the basic properties of multiple stochastic integrals.

\vskip0.5cm

\textbf{Aknowledgement} This work has been partially supported by
project by the Laboratory ANESTOC PBCT ACT-13. Karine Bertin has
been partially supported by Project FONDECYT 1090285. Soledad
Torres has been partially supported by Project FONDECYT 1070919.

\end{document}